\documentclass[reqno]{amsart}

\usepackage{amsmath}
\usepackage{amssymb}
\usepackage{graphicx}
\usepackage{xcolor}
\usepackage{amsthm}

\usepackage{tkz-tab}
\usepackage{caption}
\usepackage{latexsym}
\usepackage{subcaption} 
\usepackage{tikz}

\usepackage[colorlinks,linkcolor=blue,citecolor=blue,pdfstartview=FitH]{hyperref}

\theoremstyle{plain}
\newtheorem{thm}{Theorem}[section]

\newtheorem{lemma}[thm]{Lemma}

\theoremstyle{remark}
\newtheorem{mcase}{Motivating Case}
\newtheorem{rem}[thm]{Remark}

\newcommand{\co}{\colon\thinspace}
\newcommand{\Q}{\mathbb{Q}}
\newcommand{\R}{\mathbb{R}}
\newcommand{\Z}{\mathbb{Z}}
\newcommand{\calG}{\mathcal{G}}

\begin{document}

\title{Continued fractions and lines across the Stern--Brocot diagram}

\author{H. Abramson \and E. Chesebro \and V. Cummins \and C. Emlen \and R. Grady \and  K. Ke}

\address{Department of Mathematical Sciences, University of Montana} 
\email{Eric.Chesebro@mso.umt.edu} 

\address{Department of Mathematical Sciences, Montana State University} 
\email{ryan.grady1@montana.edu} 

\begin{abstract}   
This paper concerns the relationships between continued fractions and the geometry of the Stern-Brocot diagram.  Each rational number can be expressed as a continued fraction $[a_0; a_1, \ldots, a_n]$ whose terms $a_i$ are integers and are positive if $i \geq 1$.  Select an index $i \in \{ 1, \ldots, n \}$ and replace $a_i$ with an integer $m$ to obtain a continued fraction expansion for an extended rational $\alpha_m \in \Q \cup \{ \infty \}$.  This paper shows that the vertices of the Stern-Brocot diagram corresponding to the numbers $\{ \alpha_m \}_{m \in \Z}$ lie on a pair of (extended) Euclidean lines across the diagram.  The slopes of these two lines differ only by a sign change and they meet at the point $L=\left([a_0; a_1, \ldots, a_{i-1}], 0\right) \in \R^2$.  Moreover, as $\lvert m \rvert \to \infty$, the associated vertices move down these lines and converge to $L$.  This paper concludes with a discussion which interprets this result in the context of 2-bridge link complements and Thurston's work on hyperbolic Dehn surgery.
\end{abstract}

\maketitle

\section{Introduction}

A sequence of integers $(a_0, \ldots, a_n)$ determines $[a_0; a_1 , \ldots, a_n] \in \Q \cup \left\{ \frac{1}{0} \right\}$ by the continued fraction
\[ a_0+\cfrac{1}{a_1+\cfrac{1}{\begin{array}{ccc} a_2 + & &  \\ & \ddots & \\ & & +\frac{1}{a_n}. \end{array}}}. \]
Its $j^\text{th}$ {\em convergent} is the number $[a_0; a_1, \ldots, a_j]$.  There is also a geometric graph $\calG \subset \R^2$ called the {\em Stern-Brocot diagram}, whose vertices are parametrized by $\Q$.  This paper describes certain connections between continued fractions, convergents, and the geometry of $\calG$.

\begin{figure}
   \centering
   \includegraphics[width=2.75in]{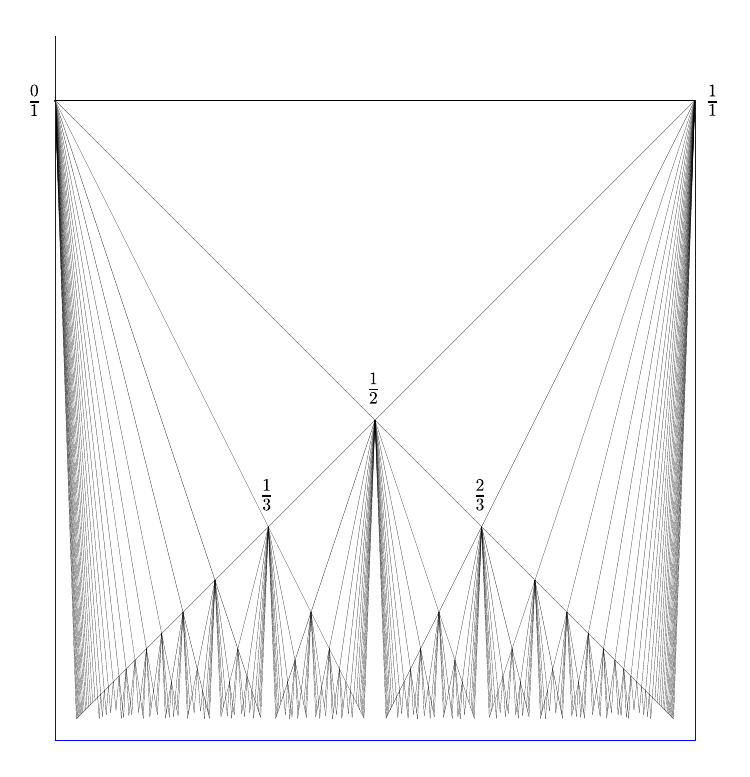} 
   \caption{A portion of the Stern--Brocot diagram $\calG$.}   \label{fig: G}
\end{figure}

This paper will take the convention that elements of $\Q-\Z$ are written in lowest terms with positive denominators.  Sometimes integers will be written as quotients with unit denominator.  A pair of rationals $\{ p/q, r/s \}$ is called a {\em Farey pair} if $ps-rq =\pm 1$.  A {\em Farey triple} is a triple of rational numbers which are pairwise Farey pairs.  

Let $\nu \co \Q \to \R^2$ be the function $\nu(p/q)=\left( p/q, 1/q\right)$.  The vertices of the Stern-Brocot diagram $\calG$ are defined as the points $\nu(\Q)$.  Two vertices of $\calG$ are connected by an edge if and only if their preimages under $\nu$ form a Farey pair.  In \cite{Hat}, Hatcher shows that the interiors of the set of all edges in $\calG$ are mutually disjoint.  This is evident in Figure \ref{fig: G}, which depicts a portion of $\calG$ over the vertices $\nu \left( \Q \cap [0,1] \right)$.  Just as Farey pairs can be viewed as edges in $\calG$, every Farey triple appears as a triangle in $\calG$.  Since every Farey pair belongs to exactly two Farey triples, every edge of $\calG$ is shared by exactly two triangles.  

\begin{figure}
   \centering
   \includegraphics[width=4.95in]{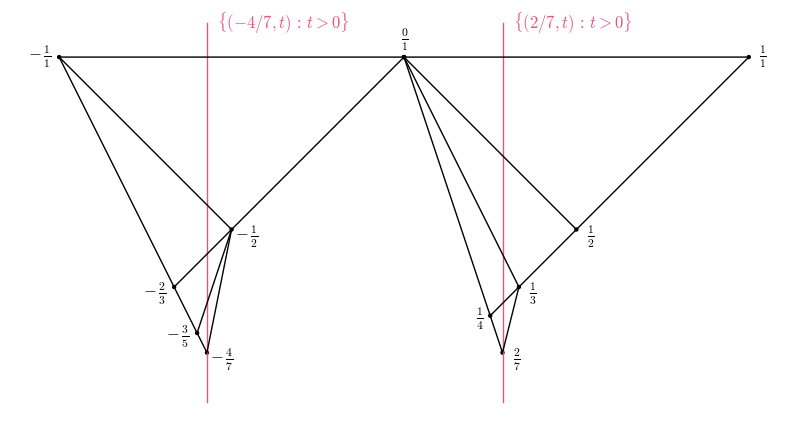} 
   \caption{The funnels $F_{-4/7}$ and $F_{2/7}$ correspond to the continued fraction expansions $-4/7=[-1;2,3]$ and $2/7=[0;3,2]$.}   \label{fig: fun}
\end{figure}

One connection between continued fractions and the geometry of $\calG$ can made through the funnels of $\calG$.  Given $\alpha \in \Q$, the {\em funnel} $F_\alpha$ is the union of all triangles in $\calG$ which meet the the ray $\{ (\alpha, t) \, : \, t \in \R^+\}$.  For instance, the funnels $F_{-4/7}$ and $F_{2/7}$ are shown in Figure \ref{fig: fun}.  If $\omega$ is a vertex of $F_\alpha$, its {\em index} is the number of edges of $F_\alpha$ which have $\omega$ as an endpoint and meet the defining ray for $F_\alpha$.  

For $\alpha \in \Q$, the Euclidean algorithm guarantees a sequence $(a_0, \ldots, a_n)$ of integers with $\alpha = [a_0; a_1,\ldots, a_n]$, $a_j\geq 1$ whenever $j\geq 1$, and $a_n \geq 2$.  When a sequence $(a_0, \ldots, a_n)$ has these properties, it is called {\em standard}.  The following theorem relates funnels and convergents.   

\begin{thm} [Theorem 2.1 of \cite{Hat}] \label{thm: hat}
Let $c_j$ be the $j^\text{th}$ convergent for a standard sequence of integers $(a_0, \ldots, a_n)$.
\begin{enumerate}
\item If $j$ is even, then $\nu(c_j)$ lies on the left edge of $F_\alpha$.  Otherwise, $\nu(c_j)$ lies on the right edge.
\item The index of $\nu(c_0)$ is $a_1$ and the index of $\nu(c_{n-1})$ is $a_n$.
\item If $0<j<n-1$, then the index of $\nu(c_j)$ is $1+a_j$.
\end{enumerate}
\end{thm}

This theorem shows that the continued fraction expansion for $\alpha$ can be determined immediately by looking at a picture of its  funnel.  Conversely, the combinatorial properties of the funnel can be deduced from the continued fraction expansion.

Suppose $(a_0, \ldots, a_n) \in \Z^{n+1}$ is standard and fix $i \in \{ 1, \ldots, n\}$.  Then every $m\in \Z$ determines $\alpha_m \in \Q \cup \{ \frac{1}{0} \}$ by 
\[ \alpha_m = [a_0; a_1 ,\ldots, a_{i-1}, m, a_{i+1}, \ldots, a_n].\]

\begin{figure}
   \centering
   \includegraphics[width=4.75in]{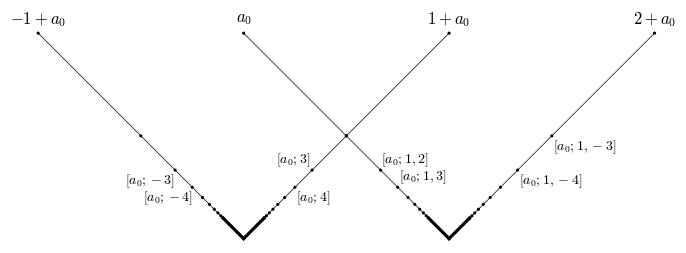} 
   \caption{Motivating cases \ref{mcase 2} and \ref{mcase 3} concern the numbers $[a_0; m]$ and $[a_0; 1, m]$.  The corresponding Stern-Brocot vertices appear along Euclidean line segments as shown above.   } \label{fig: fun} 
\end{figure}

\begin{figure}
   \centering
   \includegraphics[width=1.5in]{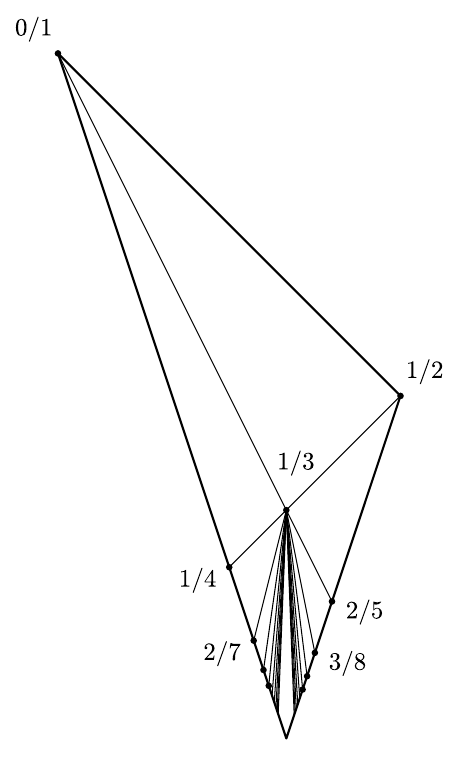} 
   \caption{The vertices in $\calG$ for $[0;3,m]=\frac{m}{1+3m}$ lie on the left edge of the triangle when $m\geq 0$.  The vertices for $[0;2,1,m]=\frac{1+m}{2+3m}$ lie on the right edge when $m\geq 0$.}   \label{fig: tri}
\end{figure}

\begin{mcase} \label{mcase 1} Suppose that $i=n$ and let $\gamma=[a_0; a_1, \ldots, a_{n-1}]$.  It follows from Theorem \ref{thm: hat} that, if $m$ is positive, then $\gamma$ and $\alpha_m$ make a Farey pair.    Moreover, using the construction of $\calG$ given in \cite{Hat}, this implies that the points $\nu(\alpha_m)$ all lie on a single Euclidean line which passes through through $(\gamma,0)$.  This can be seen in Figure \ref{fig: tri}.  \end{mcase}

\begin{mcase} \label{mcase 2}  Suppose that $n=1$ and $\gamma=a_0$.  It must be true that $i=1$ and so $\alpha_m=[a_0; m] = a_0+1/m$.  Here, if $m<0$, then $\nu(\alpha_m)$ lies on the line segment between $(-1+\gamma,1)$ and $(\gamma,0)$ and, if $m>0$, then $\nu(\alpha_m)$ lies on the line segment between $(1+\gamma,1)$ and $(\gamma,0)$.  The value $\alpha_0$ is defined as $1/0$ and $\nu(1/0)$ is defined as the added point in $\R^2 \cup \{ \infty \}$.  Hence, $\nu(\alpha_0)$ lies on the {\em extended } Euclidean lines through $(\pm 1+\gamma,1)$ and $(\gamma,0)$.  \end{mcase}

\begin{mcase} \label{mcase 3}  Suppose that $i=n=2$ and $a_1=1$.   Put $\gamma=a_0+1$ and $\alpha_m=[a_0;1,m]=a_0+m/(m+1)$.  If $m<-1$, then $\nu(\alpha_m)$ lies on the line segment between $(1+\gamma,1)$ and $(\gamma,0)$ and, if $m>-1$, then $\nu(\alpha_m)$ lies on the line segment between $(-1+\gamma,1)$ and $(\gamma,0)$.  The infinite point occurs at $m=-1$ where $\nu(\alpha_{-1})$ lies on the extended Euclidean lines as before.  \end{mcase}

The main result of this paper concerns the general case.  

\begin{thm} \label{thm: main}
Suppose $(a_0, \ldots, a_n) \in \Z^{n+1}$ is standard and fix $i \in \{ 1, \ldots, n\}$.  For each $m\in \Z$ define  
\[ \alpha_m = [a_0; a_1 , \ldots, a_{i-1}, m, a_{i+1}, \ldots, a_n] \in \Q \cup \left\{ 1/0 \right\}\]
and $\gamma = [a_0; a_1, \ldots, a_{i-1}]$.  Let $\ell^+$ be the extended Euclidean line in $\R^2 \cup \{ \infty \}$ through $(\gamma, 0)$ and $\nu(\alpha_1)$.  Let $\ell^-$ be the extended Euclidean line obtained from $\ell^+$ by reflection across the $x$-axis.  Then, for every $m\in \Z$, 
\[ \nu(\alpha_m) \in \ell^+ \cup \ell^-.\]
\begin{enumerate}
\item \label{state 1} If $\nu(\alpha_m) = \infty$, then $m \in \{ -1, 0 \}$ and, if $\nu(\alpha_0) = \infty$, then $n=1$.
\item \label{state 2} If $m \leq -2$ then $\nu(\alpha_m) \in \ell^-$ and if $m \geq 0$ then $\nu(\alpha_m) \in \ell^+$.
\item \label{state 3} $\left| \nu(\alpha_m) - (\gamma, 0) \right|$ and $\left| \nu(\alpha_{-m}) - (\gamma, 0) \right|$ converge monotonically to zero as $m \to \infty$.
\end{enumerate}
\end{thm}

\begin{rem}
Motivating Cases \ref{mcase 2} and \ref{mcase 3} show that both situations for (\ref{state 1}) of Theorem \ref{thm: main} do occur.  For another example with $\alpha_{-1}=1/0$, consider the expansion $[0;2,-1,2]$.  Also, $\nu (\alpha_{-1})$ can lie on the finite part of either $\ell^+$ or $\ell^-$.  In most of these cases, $\nu(\alpha_{-1}) \in \ell^-$, but it is not difficult to realize the other possibilities; for instance $\nu([0;2,1,-1,2])$ is in  $\ell^+$.
\end{rem}

\begin{figure}[h] 
   \centering
   \includegraphics[width=4.75in]{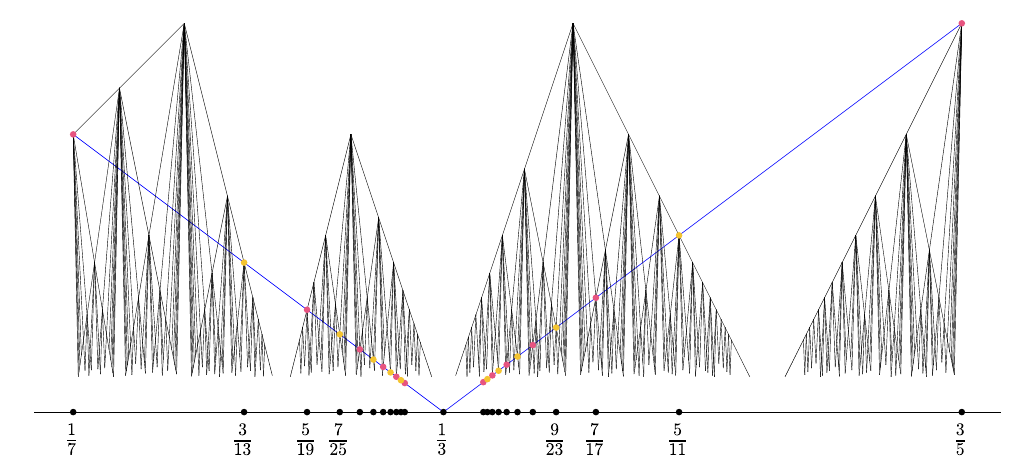} 
   \caption{This shows a portion of $\calG$ along with the (blue) lines $\ell^\pm$ for $\alpha_m=[0;3,m,4]$. The vertices for $\nu(\alpha_m)$ are shown as pink dots.  In this case, the lines $\ell^\pm$ are shared with those for the numbers $\beta_m=[0;2,1,m,4]$.  The vertices for $\nu(\beta_m)$ are shown as yellow dots.}   \label{fig: 3m4}
\end{figure} 

Our interest in this theorem came from connections to the study of the topology of knots and links in the 3-sphere.  It is well-known that continued fractions and funnels are valuable concepts when working with an important class of links called {\em 2-bridge links}. Indeed, this perspective is essential for the arguments in \cite{HT}, \cite{SW}.  We recall these connections in the final section, where we also interpret the preceding theorem in this setting.

\subsection*{Acknowledgement} 
This project was made possible by the generous support of the Center for Undergraduate Research in Mathematics and their NSF grant DMS-1722563.  RG was is supported by the Simons Foundation under Travel Support/Collaboration 9966728 during the preparation of the manuscript. 

\section{Lemmas} \label{sec: lem}

Hatcher's approach to continued fractions and the proof of Theorem \ref{thm: hat} utilizes the action of $2\times 2$ non-singular matrices with integer entries both on $\Z^2$ by ordinary matrix multiplication, as well as on $\Q \cup \left\{\frac{1}{0}\right\}$ by M{\"obius} transformations.

For a non-singular matrix $X=\left( \begin{smallmatrix} a&b\\c&d \end{smallmatrix} \right)$ with integer entries define the M{\"obius} transformation $\phi_X \co \Q\cup \left \{ \frac{1}{0} \right \} \to \Q\cup \left \{ \frac{1}{0}\right \}$ by 
\[ \phi_X(v) \ =\ \begin{cases} 
\frac{a}{c} & \text{if } v=\frac{0}{1}  \\
\frac{av+b}{cv+d} & \text{otherwise.} \end{cases}  \]
Let $\theta \co \Z^2-\{\bar{0}\} \to \Q \cup \left\{\frac{1}{0}\right\}$ be the quotient function $\theta \left( \begin{smallmatrix} p \\ q \end{smallmatrix} \right) = \frac{p}{q}$.  A quick calculation shows that 
\begin{equation} \label{eq: commute} \phi_X \theta = \theta X.\end{equation}
For $a \in \Z$, define $X_a=\left( \begin{smallmatrix} 0&1\\1&a \end{smallmatrix}\right)$.  Equation (\ref{eq: commute}) implies $\theta \left( X_a \left( \begin{smallmatrix} p \\ q \end{smallmatrix} \right) \right)  =   \left[0; a, p/q \right]$.  This simple observation is used by Hatcher in \cite{Hat} to establish the following fundamental fact which he uses to prove Theorem \ref{thm: hat}.

\begin{lemma} \label{lem: hat}
Suppose that $(a_0, \ldots, a_n)$ is a sequence of integers and 
\[ \begin{pmatrix} p \\ q \end{pmatrix} \ = \ \begin{pmatrix} 1&a_0 \\ 0&1 \end{pmatrix} \begin{pmatrix} 0&1\\1&a_1 \end{pmatrix} \cdots \begin{pmatrix} 0&1\\1&a_n \end{pmatrix}\begin{pmatrix} 0\\1 \end{pmatrix}\]
then $[a_0; \, a_1, \ldots, a_n] = p/q$.
\end{lemma}

The next lemma is essentially a restatement of Lemma \ref{lem: hat} in the case $a_0=0$.  A few extra details are singled out, as they will be handy later.  Let $I$  denote the $2\times 2$ identity matrix.  

\begin{lemma} \label{lem: properties}
Suppose $a_1, \ldots, a_j \in \Z$ and let $B_j = I X_{a_1} \cdots X_{a_j}$.  (Note that $B_0$ is defined as the identity.)
\begin{enumerate}
\item The entries of $B_j$ are integers,
\item The first column of $B_j$ is the second column of $B_{j-1}$,
\item If $B_j = \left( \begin{smallmatrix} a&c\\b&d \end{smallmatrix} \right)$, then 
$[0; a_1, \ldots , a_j]=c/d$,
\item $\det(B_j) = (-1)^j$.
\item If $\bar{v}$ is a row or column of $B_j$, then the entries of $\bar{v}$ are relatively prime unless one entry is zero.  In this case, the other entry of $\bar{v}$ must be $\pm 1$.
\item If $a_1, \dots, a_j$ are positive and $j\geq 1$, then the top entry of the second column of $B_j$ is at least the sum of the entries of the top row of $B_{j-1}$ and, similarly, the bottom entry of the second column of $B_j$ is at least the sum of the entries of the bottom row of $B_{j-1}$.
\item If $a_1, \dots, a_j$ are positive and $B_j$ has zero as an entry, then $j \in \{ 0,1 \}$.
\end{enumerate}
\end{lemma}

\proof  Properties (1) and (2) are immediate consequences of matrix multiplication. Property (3) is Lemma \ref{lem: hat} when $a_0=0$.  The determinant assertion follows because $\det(X_a)=-1$ for every number $a$. 

For (5), if one of the entries of $\bar{v}$ is zero, then (1) and (4) imply that the other entry must be $\pm 1$.  Now suppose that neither entry is zero, write $B_j$ as in (3), and assume that $\bar{v}$ is the first column of $B_j$.  Suppose also that $n$ is an integer that divides both $a$ and $b$.  Then $n$ must also divide $ad$ and $bc$.  Therefore, $n$ divides $\det(B_j)$.  Since this is $\pm 1$, so must be $n$.  The same argument works if $\bar{v}$ is the second column or a row.

To see property (6), write $B_{j-1} = \left( \begin{smallmatrix} w&x\\y&z \end{smallmatrix}\right)$.  Then
\[ B_j \ = \ \begin{pmatrix} w&x\\y&z \end{pmatrix} \begin{pmatrix} 0&1 \\ 1&a_j \end{pmatrix} \ = \ \begin{pmatrix} x&w+xa_j \\ z& y+za_j \end{pmatrix}.\]
Since $x,z\geq 0$ and $a_j>0$,
\[ w+xa_j \geq w+x \quad \text{and} \quad y+za_j \geq y+z.\]

Lastly, if $j\geq 2$, properties (2) and (6) imply that all entries of $B_2$ are at least $1$.  Property (6), in turn, implies that the entries of $B_j$ are also at least $1$.  This proves property (7).

\endproof

The next lemma clarifies when points $\nu(p/q)$ lie on certain lines.

\begin{lemma} \label{lem: det}
Suppose that $\{ a, b \}$ and $\{ p, q \}$ are relatively prime pairs of integers.  Suppose also that $t,u \in \Z$ and $bqu \neq 0$.  Then the point $\nu(a/b)$ lies on the line through $\nu(p/q)$ and $\left(t/u,0\right)$ if and only if 
\[ \begin{vmatrix} p&t\\q&u \end{vmatrix} = \begin{vmatrix} a&t\\b&u \end{vmatrix}.\]
\end{lemma}

\proof
The slope of the line through $\nu(a/b)$ and $\left(t/u,0\right)$ is
\[ \frac{1/b - 0}{a/b-t/u} = \frac{u}{au-bt}.\]
Similarly, the slope of the line through $\nu(p/q)$ and $\left(t/u,0\right)$ is $u/(pu-qt)$.  These two slopes are the same if and only if
\[ \begin{vmatrix} p&t\\q&u \end{vmatrix} = \begin{vmatrix} a&t\\b&u \end{vmatrix}.\]
\endproof

Since the proof of the main theorem will reduce to the case of expansions of rationals in $\Q \cap [0,1/2]$, the next lemma will be helpful.

\begin{lemma} \label{lem: cf}
Suppose $(a_1, \ldots a_j)$ is a sequence of positive integers.  If $a_1=1$ then
\[ \frac{1}{2} \leq [0; a_1, \ldots, a_j] \leq 1,\]
where the upper inequality is an equality if and only if $a_j=j=1$.
If $a_1 \geq 2$, then
\[ 0 < [0; a_1, \ldots, a_j] \leq \frac{1}{2}.\]
\end{lemma}

\proof
Since $[0; a_2, \ldots, a_j] \geq 0$,
\begin{equation} \label{eq: cf} [0;a_1, \ldots, a_j]=\frac{1}{a_1 + [0;a_2, \ldots a_j]} \leq \frac{1}{a_1}.\end{equation}
So, if $a_1\geq 2$, the second claim of the lemma holds.  Suppose then, that $a_1=1$.  Since $[0; a_2, \ldots a_j]$ and $[0;a_1, \ldots, a_j]$ are at most one, 
\[ \frac{1}{2} \leq \frac{1}{1+[0;a_2, \ldots, a_j]} = [0; a_1, \ldots, a_j] \leq 1.\]
Finally, if $[0;a_1, \ldots, a_j]=1$, statement (6) of Lemma \ref{lem: properties} implies that $j=1$.  Since $1/a_1=[0;a_1]=1$, it must also be true that $a_1=1$.
\endproof

\section{Proof of the main theorem}

\proof Suppose $(a_0, \ldots, a_n) \in \Z^{n+1}$ is standard and fix an integer $i \in \{ 1, \ldots, n\}$.  Define $\alpha_m = [a_0; a_1 , \ldots, a_{i-1}, m, a_{i+1}, \ldots, a_n]$ and $\gamma = [a_0; a_1, \ldots, a_{i-1}]$.  Let $\ell^+$  be the extended Euclidean line through $(\gamma, 0)$ and $\nu(\alpha_1)$ and  let $\ell^-$  be the extended Euclidean line obtained from $\ell^+$ by reflection across the $x$-axis.  

If $n=1$, the theorem follows from the discussion of Motivating Case \ref{mcase 2}.  So, assume $n \geq 2$.

If $a_0$ is adjusted by adding an integer, then the corresponding points $\nu(\alpha_m)$ and $\gamma$ are shifted by a horizontal translation by the same integer.  This means that it is enough to prove the theorem under the assumption that $a_0=0$.  (By Lemma \ref{lem: cf}, this is equivalent to assuming that $0 \leq \alpha_{a_i} \leq 1$.)

Define non-negative integers $r,s,t,u,v$, and $w$ by 
\[ \begin{pmatrix} r&t \\ s&u \end{pmatrix} \ = \ I X_{a_1} \cdots X_{a_{i-1}} \quad \text{and} \quad \begin{pmatrix} v \\ w \end{pmatrix} \ = \ X_{a_{i+1}} \cdots X_{a_n}  \begin{pmatrix} 0\\1 \end{pmatrix}.\]
Define also
\[P(m)=twm+rw+tv \quad \text{and} \quad Q(m)=uwm+sw+uv.\]
For any $m \in \Z$, matrix multiplication gives
\begin{align*}
 X_{a_1} \cdots X_{a_{i-1}} X_m X_{a_{a+1}} \cdots X_{a_n} \begin{pmatrix} 0\\1 \end{pmatrix} &= \begin{pmatrix} r&t\\s&u \end{pmatrix} \begin{pmatrix} 0&1 \\ 1&m \end{pmatrix} \begin{pmatrix} v\\w \end{pmatrix} \\
 &= \begin{pmatrix} P(m) \\ Q(m) \end{pmatrix}.
\end{align*}
By Lemma \ref{lem: hat}, $\alpha_m = P(m)/Q(m)$.  Statement (5) of Lemma \ref{lem: properties} shows that either $\alpha_m \in \{ 0/1, 1/0 \}$ or $P(m)$ and $Q(m)$ are relatively prime.

\medskip \noindent {\em Claim I. For every $m \in \Z$, 
\[ \left( \frac{P(m)}{Q(m)}, \frac{1}{Q(m)} \right) \ \in \ \ell^+.\]} \proof
Define
\[ w_m \ =\ \left( \frac{P(m)}{Q(m)}, \frac{1}{Q(m)} \right).\]
If $Q(m)=0$, then $w_m$ is the infinite point of the extended plane which is, by definition, an element of $\ell^+$.  So, assume $Q(m) \neq 0$.
Since 
\[ \begin{pmatrix} P(m) & t \\ Q(m) & u \end{pmatrix} \ =\ I X_{a_1} \cdots X_{a_{i-1}} X_m \begin{pmatrix} v&1\\w&0 \end{pmatrix}\]
and $\det\left( I X_{a_1} \cdots X_{a_{i-1}} X_m \right) = (-1)^i$, the determinant $\left| \begin{smallmatrix} P(m) & t \\ Q(m) & u \end{smallmatrix} \right|$ is independent of $m$.  Lemma \ref{lem: det} completes the proof of the claim.
\endproof

Notice that $w_m = \nu(\alpha_m)$ if and only if $Q(m) \geq 0$.  On the other hand, if $Q(m)<0$ then $\nu(\alpha_m)=\tau(w_m)$, where $\tau \co \R^2 \to \R^2$ is reflection across the $x$-axis.   This means that $\nu(\alpha_m) \in \ell^+ \cup \ell^-$ for every $m\in \Z$ and $\nu(\alpha_m) \in \ell^-$ if and only if $Q(m) \leq 0$.  The remaining task is now only to prove the three numbered statements of the theorem.

\medskip \noindent {\em Claim II. $Q$ is a strictly increasing function.}  \proof $Q$ is a linear function with slope $uw$, so it suffices to show that $uw >0$.  Certainly, neither $u$ nor $w$ can be negative (see, for instance, property (6) of Lemma \ref{lem: properties}).  Moreover, property (7) of Lemma \ref{lem: properties} implies that $uw \neq 0$. Indeed, if $u=0$ then (7) gives $i=1$ or $i=2$.  But $i=1$ means $\left( \begin{smallmatrix} r&t\\s&u \end{smallmatrix} \right)=I$ and then $u=1$.  If $i=2$ then $\left( \begin{smallmatrix} r&t\\s&u \end{smallmatrix} \right)=X_{a_1}$ and $u=a_1$.  This same idea works to show $w\neq 0$. Hence, $uw >0$ as claimed.  \endproof

The $y$-coordinate for $\nu(\alpha_m)$ is $1/|Q(m)|$.  Since $Q$ is an increasing linear function, statement (\ref{state 3}) of the theorem holds.  

Using that $Q$ is increasing, there is exactly one real number $x$ with $Q(x)=0$.  To complete the proof of the theorem, it suffices to show that $x$ lies in the open interval $(-2,0)$, since this would imply that $Q$ is positive on $[0, \infty)$ and negative on $(-\infty, -2]$.

\medskip \noindent {\em Claim III. If $i=1$ then $-1<x<0$.  If $a_1 \geq 2$, then $-2<x<0$.} \proof $Q(x) = 0$ and $uw>0$ so 
\[ -x \ =\ \frac{sw+uv}{uw} \ =\ \frac{s}{u} + \frac{v}{w}.\]
If $x=0$ then, since $u$ and $w$ are positive, $s$ and $v$ must both be zero.  But then, property (7) of Lemma \ref{lem: properties}, implies that $i=1=n$, which is not the situation.  Hence $x<0$.  By Lemma \ref{lem: hat}, $v/w=[0; a_{i+1}, \ldots, a_n]$ and Lemma \ref{lem: cf} gives that $0 \leq v/w \leq 1$.

If $i=1$, then $s/u=0$ and $x=-v/w$.  Since $n > 1$ and $a_n \geq 2$, Lemma \ref{lem: cf} gives \[-1<-v/w=x<0.\]

If $i>2$, then properties (6) and (7) of Lemma \ref{lem: properties} imply that $u$ is the sum of $s$ with a positive number and so $s/u<1$.  If $i=2$, then $s/u=1/a_1 \leq 1/2$.  Therefore, if $i \geq 2$, then $0 \leq s/u < 1$ and $-2 < x$.
\endproof

By Claim III, the theorem is proven unless $i \geq 2$ and $a_1=1$.  Suppose this is the case and let $\phi = \phi_X$, where $X = \left( \begin{smallmatrix} -1&1\\0&1 \end{smallmatrix} \right)$.  Define
\[ \beta_m =  \begin{cases} [0; 1+a_2, a_3, \ldots, m, \ldots, a_n] & \text{if } i>2 \\ 
[0; m, a_3, \ldots, a_n]  & \text{if } i=2. \end{cases}.\]
Because
\[ \begin{pmatrix} -1&1\\0&1 \end{pmatrix} \begin{pmatrix} 0&1\\1&a_1 \end{pmatrix} \ = \ \begin{pmatrix} 0&1\\1&1 \end{pmatrix} \begin{pmatrix} 0&1\\1&-1+a_1 \end{pmatrix},\]
the discussion of M{\"obius} transformations at the start of Section \ref{sec: lem} implies that
\[ \phi \left( \alpha_m \right) \ = \ [0;1,-1+a_1, a_2, \ldots,m, \ldots, a_n].\]
Since $a_1=1$,
\begin{align*}
\phi \left( \alpha_m \right) &=  [0;1,0, a_2, \ldots, m, \ldots, a_n]\\
&= [0; 1+a_2, a_3, \ldots, a_n].
\end{align*}
So, if $i>2$, then $\phi(\alpha_m)=\beta_m$ and, if $i=2$, then $\phi(\alpha_m)=\beta_{1+m}$.
Since $\phi$ preserves denominators, Claim III applies in the first case to show that $-2<x<0$.  In the second case, Claim III gives $0<x<-1$ as needed.
\endproof

\section{Discussion: Continued Fractions and 2-Bridge Links}

\begin{figure}[h] 
   \centering
   \includegraphics[width=3in]{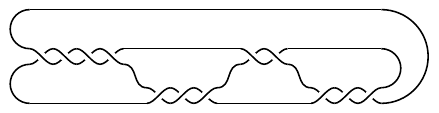} 
   \caption{The 4-plat diagram $D(4,3,2,3)$ is standard. }   \label{fig: SW}
\end{figure}

A sequence $(a_1, \ldots, a_n) \in \Z^n$ determines a 4-plat link diagram $D=D(a_1, \ldots, a_n)$ as shown in Figure \ref{fig: SW}.  The diagram is primarily made up of  two rows of horizontal twist regions, each of which consists of a chain of bigons.  The left most twist region lies in the top row and has $|a_1|$ crossings.  If $a_1$ is positive, the crossings twist in the manner of a right handed screw.  Otherwise, the crossings in this first twist region are left handed.  The next twist region corresponds to $a_2$ and lies in the bottom row.  It has $|a_2|$ right handed twists if $a_2$ is negative and left handed twists if $a_2$ is positive.  In general, the $j^\text{th}$ twist region lies in the top row exactly when $j$ is odd.  Moreover, top row twists are right handed when $a_j$ is positive and bottom row twists are right handed when $a_j$ is negative.  Twist regions are connected into a link diagram as indicated in Figure \ref{fig: SW}.  If each $a_j$ is positive and $a_1, a_n \geq 2$, then $D$ is called {\em standard}.  Note that standard diagrams are alternating in the sense that, while following the strands of the link, crossings alternate over and under.

As is common, this paper interprets a link diagram as describing an embedding of a collection of circles into $\R^3$ (up to isotopy).  A pair of link diagrams are considered to be equivalent if the differ by a finite sequence of Reidemeister moves (see \cite{A}) and/or a simultaneous reflection of all crossings in the diagram.  For instance, the diagrams $D(a_1, \ldots, a_n)$ and $D(-a_1, \ldots, -a_n)$ are equivalent.  The equivalence class of a link diagram is called a link and a link represented by a 4-plat diagram is called a {\em 2-bridge link}.

Suppose that $(a_1, \ldots, a_n) \in \Z^n$ and $(b_1, \ldots, b_m) \in \Z^m$ and let $p/q=[0;a_1, \ldots, a_n]$ and $p'/q'=[0;b_1, \ldots, b_m]$.  Schubert showed in \cite{Sch} that the 2-bridge link diagrams $D(a_1, \dots, a_n)$ and $D(b_1, \ldots b_m)$ are equivalent if and only if $q'=q$ and $p'$ is equivalent modulo $q$ to either $p$, $-p$, $p^{-1}$ or $-p^{-1}$.  See also \cite{BZ}, \cite{Con}, and \cite{KL} for more on this result.  It follows that every 2-bridge link is represented by a standard 4-plat diagram.  Moreover, every rational number $p/q$ in the interval $[0,1/2]$ determines a well-defined 2-bridge link $L(p/q)$ and every 2-bridge link arises in this manner. 

Consider a standard link diagram $D(a_1, \ldots, a_n)$.  Then $(0,a_1, \ldots, a_n)$ is a standard sequence of integers.  As usual, fix $i \in \{ 1, \ldots, n\}$ and set $\alpha_m = [0; a_1, \ldots, m, \ldots, a_n]$.  The diagram $D(a_1, \ldots,m, \ldots, a_n)$ is a diagram for the 2-bridge link $L(\alpha_m)$.  

This setup is of interest for topologists.  Very often knots and links are studied in terms of the topology of their {\em complements}.  Recall that the link diagram determines an (isotopy class) of embedded circles in $\R^3$.  In turn, $\R^3$ sits naturally as a subset of the 3-sphere $S^3=\R^3 \cup \{ \infty\}$.   The complement of the link is the topological space obtained by deleting the embedded circles from $S^3$.  Standard 3-manifold topology shows that the topology of the link complement is independent of the diagram chosen to represent the link.  This creates an intimate connection between the study of links and that of 3-dimensional spaces.  This connection is particularly strong for the subclass of knot complements, where the link embedding is of a single loop.  In this setting, Gordon and Lueke  famously showed that two knot diagrams are equivalent if and only if their complements have the same topology (see \cite{GL1} and \cite{GL2}).

This setup is also interesting from the point of view of geometry.  If $n \geq 2$ and $|m|$ is not too small then, using Thurston's work \cite{Th_notes}, Menasco showed in \cite{Me} that the complement of $L(\alpha_m)$ admits a complete hyperbolic metric under which the volume of complement is finite.  This metric is uniquely determined by the topology of the complement.  Furthermore, Thurston's hyperbolic Dehn surgery theorem implies that these geometric spaces are intimately related.  In fact, as $|m| \to \infty$ the hyperbolic complements of the links $L(\alpha_m)$ converge to a fixed hyperbolic link complement.  

The geometries of the hyperbolic 2-bridge link complements are succinctly described in \cite{C2} in terms of the positions of their corresponding vertices $\nu(\alpha_m)$.  Theorem \ref{thm: main} says that, not only do the points $\nu(\alpha_m)$ determine the geometry of the corresponding hyperbolic links but also, as these link complements converge their Stern-Brocot vertices also converge down the lines $\ell^\pm$.

\bibliographystyle{plain}
\bibliography{MF}
\end{document}